# METHODOLOGY AND CONVERGENCE RATES FOR FUNCTIONAL LINEAR REGRESSION

By Peter Hall[1] and Joel L. Horowitz[2]

*Australian National University and Northwestern University*

In functional linear regression, the slope "parameter" is a function. Therefore, in a nonparametric context, it is determined by an infinite number of unknowns. Its estimation involves solving an ill-posed problem and has points of contact with a range of methodologies, including statistical smoothing and deconvolution. The standard approach to estimating the slope function is based explicitly on functional principal components analysis and, consequently, on spectral decomposition in terms of eigenvalues and eigenfunctions. We discuss this approach in detail and show that in certain circumstances, optimal convergence rates are achieved by the PCA technique. An alternative approach based on quadratic regularisation is suggested and shown to have advantages from some points of view.

**1. Introduction.** In functional linear regression, data pairs $(X_i, Y_i)$ are generated by the model

$$(1.1) \qquad Y_i = a + \int_{\mathcal{I}} b X_i + \varepsilon_i, \qquad 1 \leq i \leq n.$$

The $X_i$'s are random functions, $\mathcal{I}$ denotes the interval on which each such function is defined, the intercept $a$ and the errors $\varepsilon_i$ are scalars and the slope $b$, the main object of our interest in this paper, is a function. The model (1.1) is applicable in a wide range of settings, including many where data are becoming available only through new developments in technology.

For example, in near-infrared spectroscopy applied to data on different cereal-grain types (e.g., different varieties of wheat), $X_i(t)$ denotes the in-

Received November 2004; revised March 2006.
[1]Supported in part by an ARC grant.
[2]Supported in part by NSF Grants SES-99-10925 and SES-03-52675.

tensity of reflected radiation recorded at the spectrometer when the wavelength equals $t$ and $Y_i$ denotes the level of a particular protein for the $i$th cereal type. By constructing the linear regression at (1.1), we can predict, from data on a new function $X$, the level of protein for that cereal type. This is especially useful in practice, since the explanatory variables $X_i$ are very easy and inexpensive to observe in the field using hand-held equipment, whereas direct calculation of the $Y_i$ requires expensive and time-consuming analysis in a laboratory. There is an extensive literature on this problem; see, for example, [26, 31].

Once an estimator $\hat{b}$ of the slope $b$ is available, it is straightforward to estimate the intercept $a$, for example, as the average of the values of $Y_i - \int_{\mathcal{T}} \hat{b} X_i$. Therefore, much interest in the literature focuses on estimating $b$. The conventional approach, discussed, for example, by Ramsay and Silverman ([24], Chapter 10 and [25]), is based on principal components analysis or PCA. Although this method has been widely discussed (e.g., [3, 7, 14]), relatively little is known about convergence rates of estimators, apart from upper bounds. In this paper, we shall give optimal convergence rates in this problem and discuss PCA-based estimators which attain those rates. The known upper bounds for convergence rates are an order of magnitude greater than the minimax-optimal rates derived in this paper.

An alternative approach based on Tikhonov, or quadratic, regularization [29] will also be addressed. To the best of our knowledge, this approach has not been considered before in functional data analysis, although it has been widely applied to the solution of other ill-posed problems. In particular, quadratic regularisation methods are increasingly studied in the statistics literature; see, for example, work of Efromovich and Koltchinskii [11] and Cavalier et al. [5] on optimality properties.

We shall show that the Tikhonov regularisation approach is also able to achieve optimal convergence rates and that it is robust against potential problems caused by tied, or closely spaced, eigenvalues in the spectral decomposition on which PCA is based. The difficulties that close eigenvalues can cause for PCA will be discussed using an example.

The estimation of slope and intercept parameters in functional linear regression has points in common with a range of smoothing and deconvolution problems where dimension reduction is involved; see, for example, [9, 12, 13, 28]. Work on statistical smoothing is particularly extensive and relatively well known to readers, so we shall not attempt to survey it here. The problem of estimating the slope in functional linear regression is also related to that of estimating the point-spread function in image analysis when the true image, or test pattern, is known. Here, too, significant work has been done; see, for example, [18, 32].

Of course, the literature on linear inverse problems is very much larger than this. In the statistics setting, it includes the work of Donoho [10] and Johnstone [17], who used wavelet and vaguelette methods, and that of van Rooij and Ruymgaart [30] on optimal convergence rates. There is also closely related work in economics on the subject of panel data [16], covariate measurement error [19] and estimation with instrumental variables (e.g., [2, 8, 15, 22, 23]. In statistics, there is related work on errors-in-variables problems (e.g., [4]). There is a small, but increasing, literature on applications of functional regression to longitudinal data analysis; see, for example, [6, 27].

**2. Methodology.** We shall assume that we observe independent and identically distributed data $(X_1, Y_1), \ldots, (X_n, Y_n)$, where each explanatory variable $X_i$ is a square integrable random function on the compact interval $\mathcal{I}$. The response variables $Y_i$ are generated by the model (1.1). It will be supposed that the errors $\varepsilon_i$ are independent and identically distributed with finite variance and zero mean and that the errors are also independent of the explanatory variables. Our goal is to discuss estimators of $b$ and to describe the rate at which they converge to the true function.

We begin by describing standard functional linear regression methodology, as discussed by, for example, Ramsay and Silverman ([24], Chapter 10). It is founded on spectral expansions of both the covariance of $X$ and its estimator and is constructed as follows.

Let $(X, Y, \varepsilon)$ denote a generic $(X_i, Y_i, \varepsilon_i)$ and put $K(u, v) = \mathrm{cov}\{X(u), X(v)\}$, $\bar{X} = n^{-1} \sum_i X_i$ and

$$\widehat{K}(u, v) = \frac{1}{n} \sum_{i=1}^{n} \{X_i(u) - \bar{X}(u)\}\{X_i(v) - \bar{X}(v)\}.$$

Write the spectral expansions of $K$ and $\widehat{K}$ as

$$(2.1) \quad K(u,v) = \sum_{j=1}^{\infty} \kappa_j \phi_j(u) \phi_j(v), \qquad \widehat{K}(u,v) = \sum_{j=1}^{\infty} \hat{\kappa}_j \hat{\phi}_j(u) \hat{\phi}_j(v),$$

where

$$(2.2) \qquad \kappa_1 > \kappa_2 > \cdots > 0, \qquad \hat{\kappa}_1 \geq \hat{\kappa}_2 \geq \cdots \geq 0$$

are the eigenvalue sequences of linear operators with kernels $K$ and $\widehat{K}$, respectively, and $\phi_1, \phi_2, \ldots$ and $\hat{\phi}_1, \hat{\phi}_2, \ldots$ are the respective orthonormal eigenvector (in fact, eigenfunction) sequences. We interpret $(\hat{\kappa}_j, \hat{\phi}_j)$ as an estimator of $(\kappa_j, \phi_j)$.

During the review process, it was suggested that the case where $\sum_j \kappa_j$ diverges might be explored. For example, the context $\kappa_j \sim j^{-\alpha}$, with $\alpha$ close to either 0 or $\frac{1}{2}$, might provide particular challenges. We agree that this setting is of mathematical interest. However, it should be noted that if $\mathrm{var}\, X(t)$ is bounded in $t$, then $\sum_j \kappa_j < \infty$. The case of unbounded covariance does not commonly arise in applied work.

Both sequences $\{\phi_j\}$ and $\{\hat{\phi}_j\}$ are complete in the class of square integrable functions on $\mathcal{I}$. The fact that each $\kappa_j$ is strictly positive implies that the linear operator corresponding to $K$, which takes a function $\phi$ to $K\phi$ and is defined by $(K\phi)(u) = \int K(u,v)\phi(v)\,dv$, is strictly positive definite. (To simplify notation, we use the symbol $K$ for both the kernel and the operator.) We determine the signs of $\phi_j$ and $\hat{\phi}_j$, in cases where signs are important, by insisting that $\int_\mathcal{I} \hat{\phi}_j \phi_j \geq 0$. This can be done without loss of generality, for example, by changing the sign of $\hat{\phi}_j$ to match that of $\phi_j$, since switching the signs of $\phi_j$ and $\hat{\phi}_j$ results in commensurate changes of sign for generalized Fourier coefficients such as the quantities $\hat{b}_j$ and $\hat{g}_j$ which we shall introduce below. Therefore, $\int_\mathcal{I} \hat{\phi}_j \phi_j > 0$ can be assumed without altering the values taken by estimators.

A model equivalent to (1.1) is

$$Y_i - \mu = \int_\mathcal{I} b(X_i - x) + \varepsilon_i, \qquad 1 \leq i \leq n,$$

where $x = E(X_i)$ and $\mu = E(Y_i) = a + \int bx$, with $x$ denoting a deterministic function on $\mathcal{I}$. It follows that if we define $g(u) = E[(Y - \mu)\{X(u) - x(u)\}]$, where $(X, Y)$ represents a generic pair $(X_i, Y_i)$, then

$$Kb = g.$$

Moreover, if we write $b = \sum_j b_j \phi_j$ and $g = \sum_j g_j \phi_j$, then $b_j = \kappa_j^{-1} g_j$. This suggests the estimator

$$\hat{b}(u) = \sum_{j=1}^{m} \hat{b}_j \hat{\phi}_j(u), \tag{2.3}$$

where the truncation point $m$ is a smoothing parameter, $\hat{b}_j = \hat{\kappa}_j^{-1} \hat{g}_j$, $\hat{g}_j = \int \hat{g} \hat{\phi}_j$,

$$\hat{g}(u) = \frac{1}{n} \sum_{i=1}^{n} (Y_i - \bar{Y})\{X_i(u) - \bar{X}(u)\} \tag{2.4}$$

and $\bar{Y} = n^{-1} \sum_i Y_i$.

Next, we suggest an alternative method which uses a ridge parameter $\rho$ rather than the cutoff $m$ as the smoothing parameter. Let $\widehat{K^+} = (\widehat{K} + \rho I)^{-1}$ denote the inverse of the operator $\widehat{K} + \rho I$, where $\rho > 0$ and $I$ is the identity operator. Define

$$\tilde{b} = \widehat{K^+}\hat{g} = \frac{1}{n} \sum_{i=1}^{n} (Y_i - \bar{Y})\widehat{K^+}\{X_i(u) - \bar{X}(u)\}, \tag{2.5}$$

where $\hat{g}$ is as in (2.4). Then $\tilde{b}$ is an estimator alternative to $\hat{b}$.

**3. Theoretical properties.** First, we treat the standard functional linear regression estimator $\hat{b}$, defined in (2.3). The Karhunen–Loève expansion of the random function $X$ is given by

$$X - E(X) = \sum_{j=1}^{\infty} \xi_j \phi_j,$$

where the random variables $\xi_j = \int_{\mathcal{I}} (X - EX)\phi_j$ have zero means and variances $E(\xi_j^2) = \kappa_j$ and are uncorrelated. Let $C > 1$ denote a constant. Concerning the distributions of the random function $X$ and the errors $\varepsilon$ in the model at (1.1), we shall assume that

(3.1) $X$ has finite fourth moment, in that $\int_{\mathcal{I}} E(X^4) < \infty$; $E(\xi_j^4) \leq C\kappa_j^2$ for all $j$, and the errors $\varepsilon_i$ are identically distributed with zero mean and variance not exceeding $C$.

Of the eigenvalues $\kappa_j$, we require that

(3.2) $$\kappa_j - \kappa_{j+1} \geq C^{-1} j^{-\alpha-1} \qquad \text{for } j \geq 1.$$

This condition prevents the spacings between adjacent order statistics from being too small. It also implies a lower bound on the rate at which $\kappa_j$ decreases: $\kappa_j$ must not be less than a constant multiple of $j^{-\alpha}$. The importance of (3.2) in ensuring Theorem 1, below, will be discussed following Theorem 2.

Of the Fourier coefficients $b_j$ and exponents $\alpha$ and $\beta$, we suppose that

(3.3) $$\begin{aligned}|b_j| &\leq Cj^{-\beta}, \\ \alpha > 1, \qquad &\tfrac{1}{2}\alpha + 1 < \beta.\end{aligned}$$

The first part of (3.3) can be viewed in at least two ways: as a definition of $\beta$, in terms of a given sequence $b_j$, or as a condition that the generalized Fourier coefficients $b_j$ do not decrease too quickly. The basis with respect to which these coefficients are defined is determined by the context of the problem and, more particularly, by the covariance function $K$, rather than outside the problem. This is not unnatural, for at least two related reasons. First, the basis $\phi_1, \phi_2, \ldots$ is canonical in the functional-data problem since it is the unique basis with respect to which the function $X$ can be expressed as a generalized Fourier series (its Karhunen–Loève expansion) with uncorrelated coefficients. It gives the most rapidly convergent representation of $X$ when speed of convergence is defined in an $L_2$ sense. Second, as discussed in Section 1, the representation with respect to this basis is fundamental to the most popular method for estimating $b$ and is therefore particularly deserving of study.

Note that the assumption that $K$ is bounded, or even the milder condition $\int_{\mathcal{I}} \text{var}\{X(u)\}\,du < \infty$, entails $\sum_j \kappa_j < \infty$. Further, note that (3.2) implies $\kappa_j \geq Cj^{-\alpha}$ for some constant $C > 0$. Therefore, boundedness of $K$ and (3.2) together imply that $\alpha > 1$, which is the second part of (3.3). The assumption $\frac{1}{2}\alpha + 1 < \beta$ in (3.3) requires that the function $b$ be sufficiently smooth relative to $K$, where smoothness of $K$ is expressed relative to the spectral decomposition of this function. (More concisely, $b$ should be sufficiently smooth relative to the lower bound on the smoothness of $K$ that is implied by the condition $\kappa_j \geq C j^{-\alpha}$.) Since $\alpha > 1$, a sufficient condition for $\frac{1}{2}\alpha + 1 < \beta$ is $\alpha \leq \beta$, which can be interpreted as requiring that the function $b$ be no less smooth than the lower bound on the smoothness of $K$ implied by (3.2).

Of the tuning parameter $m$, we assume that

$$m \asymp n^{1/(\alpha+2\beta)}. \tag{3.4}$$

In (3.4), the relation $r_n \asymp s_n$, for positive $r_n$ and $s_n$, means that the ratio $r_n/s_n$ is bounded away from zero and infinity.

Let $\mathcal{F}(C, \alpha, \beta)$ denote the set of distributions $F$ of $(X, Y)$ that satisfy (3.1)–(3.3) for given values of $C$, $\alpha$ and $\beta$. Let $\mathcal{B}$ denote the class of measurable functions $\bar{b}$ of the data $(X_1, Y_1), \ldots, (X_n, Y_n)$ generated by (1.1). We shall frame our next result in terms of minimax bounds. Below, the upper bound (3.5) shows performance of $\hat{b}$ and the lower bound (3.6) reflects performance of any estimator of $b$. The fact that the convergence rate is the same in each instance implies that the rate for $\hat{b}$ is optimal in a minimax sense.

THEOREM 1. *If* (3.1)–(3.4) *hold, then*

$$\lim_{D \to \infty} \limsup_{n \to \infty} \sup_{F \in \mathcal{F}} P_F \left\{ \int_\mathcal{I} (\hat{b} - b)^2 > D n^{-(2\beta-1)/(\alpha+2\beta)} \right\} = 0 \tag{3.5}$$

*as* $n \to \infty$. *Furthermore,*

$$\liminf_{n \to \infty} n^{(2\beta-1)/(\alpha+2\beta)} \inf_{\bar{b} \in \mathcal{B}} \sup_{F \in \mathcal{F}} \int_\mathcal{I} E_F (\bar{b} - b)^2 > 0. \tag{3.6}$$

It follows from (3.5) that for each $F \in \mathcal{F}$,

$$\int_\mathcal{I} (\hat{b} - b)^2 = O_p(n^{-(2\beta-1)/(\alpha+2\beta)}).$$

The theorem is proved in Section 5. The fact that (3.5) is expressed in terms of a probability rather than an expected value is not significant. By modifying the estimator $\hat{b}$ using a truncation point, to prevent $\hat{b}$ taking values that are too large, we may state and prove (3.5) in the more traditional form; compare (3.10) below. We do not do this, since the present form of $\hat{b}$ is the one actually used by statisticians.

Convergence rates of the form $n^{-(2\beta-1)/(\alpha+2\beta)}$ are generic to a large class of noisy inverse problems where the difficulty of inverting the operator is an increasing function of $\alpha$ and the smoothness of the target function is an increasing function of $\beta$. For example, this rate arises in the context of problems discussed by Cavalier et al. [5]. See equation (7) there and note that the appropriate values of the components of that formula are $\lambda_i = 1$ for $1 \leq \lambda_i \leq m$ and $\lambda_i = 0$ otherwise, $\theta_i = b_i$, $\sigma_i^2 = \text{var}(\xi_i)$ and $\varepsilon^2 = n^{-1}$. Of course, Theorem 1 cannot be derived from the results of Cavalier et al. [5], but, since the problem is of the same broad type, the rates enjoy the same form and have exactly the same formula if we make the substitutions above. Connections of this nature are frequently highlighted in the literature, for nonlinear inverse problems (see, e.g., [20, 21]) as well as linear ones. In particular, similar remarks can be made about the rates given by Hall and Horowitz [15].

Next, we address the alternative estimator $\tilde{b}$ in (2.5), where the smoothing parameter is the ridge $\rho$, rather than the cutoff $m$. Assumptions (3.2)–(3.4) are replaced by

$$j^{-\alpha} \leq C\kappa_j, \tag{3.7}$$

$$|b_j| \leq Cj^{-\beta}, \qquad \alpha > 1, \qquad \alpha - \tfrac{1}{2} < \beta, \tag{3.8}$$

$$\rho \asymp n^{-\alpha/(\alpha+2\beta)}, \tag{3.9}$$

respectively. Let $\mathcal{G}(C, \alpha, \beta)$ denote the set of distributions $F$ of $(X, Y)$ that satisfy (3.1), (3.7) and (3.8) for given values of $C$, $\alpha$ and $\beta$.

The result below is a direct analogue of Theorem 1 in the case of $\tilde{b}$ rather than $\hat{b}$, except that we replace the probability bound (3.5) by one on expected value.

THEOREM 2. *If* (3.1) *and* (3.7)–(3.9) *hold, then*

$$\sup_{F \in \mathcal{G}} \int_{\mathcal{I}} E_F(\tilde{b} - b)^2 = O(n^{-(2\beta-1)/(\alpha+2\beta)}) \tag{3.10}$$

*as* $n \to \infty$. *Furthermore,*

$$\liminf_{n \to \infty} n^{(2\beta-1)/(\alpha+2\beta)} \inf_{\bar{b} \in \mathcal{B}} \sup_{F \in \mathcal{G}} \int_{\mathcal{I}} E_F(\bar{b} - b)^2 > 0. \tag{3.11}$$

A proof of (3.10) can be developed along the lines of that of Theorem 4.1 of Hall and Horowitz [15] and so will not be given here; a proof of (3.11) is identical to that of (3.6). There is no close connection between the convergence rates in (3.10) and those in [15]. In fact, the only significant linkage is that both rates are obtained by using Tikhonov regularisation to solve a linear inverse problem. From a conventional statistical viewpoint, our work is much closer to that of linear regression in a large number of dimensions than it is to instrumental variables problems.

Condition (3.7) is weaker than (3.2). For example, the latter excludes cases where two or more of the eigenvalues $\kappa_j$ are close together, in particular, where they are tied. [When employing the approach (2.5), it is not necessary to assume strict inequality among the $\kappa_j$'s.] Indeed, if closely spaced eigenvalues are permitted, (3.5) in Theorem 1 can fail while (3.10) in Theorem 2 holds. This is perhaps best illustrated by an example, which we give below, in a setting where there are long strings of tied eigenvalues. The assumption of perfect ties can be relaxed by permitting the $\kappa_j$'s to be very close to one another, but not identical. The argument there is more complex, however.

Let $\gamma, \tau$ denote constants satisfying $1 < \gamma \leq \alpha\tau$ and let $j_k$ equal the least integer not less than $k^{k\tau}$. Put $\mathcal{J}_k = \{j_k, j_k + 1, \ldots, j_{k+1} - 1\}$ and define $\kappa_j = k^{-k\gamma}$ for all $j \in \mathcal{J}_k$. Then for $j$ in this range,

$$(3.12) \qquad \kappa_j = k^{-k\gamma} \geq k^{-k\alpha\tau} \geq j_k^{-\alpha} \geq j^{-\alpha}$$

and also, $j_{k+1}/j_k \sim e^\tau k^\tau$ as $k$ increases. Property (3.12) implies (3.7), but (3.2) fails because of the ties.

Those ties mean that the functions $\phi_j$, for $j$ in the block $\mathcal{J}_k$, are not even identifiable. Indeed, any permutation of the function sequence $\phi_j$, $j \in \mathcal{J}_k$, is equally appropriate, since within-block permutations of the $\phi_j$'s do not lead to violations of the condition that the $\kappa_j$'s are nondecreasing. For the same reason, while the (unordered) set of function estimators, $\widehat{\Phi}_k = \{\hat{\phi}_j : j \in \mathcal{J}_k\}$, converges to the set $\Phi_k = \{\phi_j : j \in \mathcal{J}_k\}$ as $n \to \infty$, for each $k$, the individual estimators $\hat{\phi}_j$ are not consistent for the respective functions $\phi_j$.

If the sum in (2.3) is taken over a whole number of blocks $\mathcal{J}_k$, this inconsistency does not cause any difficulties in estimating the slope function $b$. There are problems, however, if the integer $m$ in (2.3) falls midway through one of the blocks $\mathcal{J}_k$. For definiteness, take $m$ to equal the integer part of $n^{1/(\alpha+2\beta)}$, thereby satisfying (3.4). Define $k_0 = k_0(n)$ to be the unique value of $k$ such that $m \in \mathcal{J}_k$. Then along an infinite sequence, $\mathcal{N}$ say, of values of $n$, we have

$$(3.13) \qquad \tfrac{1}{2}(j_{k_0} + j_{k_0+1}) \leq m \leq j_{k_0+1} - 1.$$

Condition (3.13) ensures that the set of integers $j \in \mathcal{J}_{k_0}$ that lie between $j_{k_0}$ and $m$ comprises at least half of $\mathcal{J}_{k_0}$. Moreover, since $j_{k+1}/j_k \sim e^\tau k^\tau$, then for all sufficiently large $n \in \mathcal{N}$, the value of

$$\frac{1}{m}\#\{j : j \in \mathcal{J}_{k_0} \text{ such that } 1 \leq j \leq m\}$$

converges to 1 as $n \to \infty$. We shall call these properties (P).

An argument based on symmetry shows that if $p = \hat{p}$ is the random permutation of $\mathcal{J}_{k_0}$ defined to minimize any given symmetric measure of performance of $\widehat{\Phi}$ as an estimator of $\Phi_j$, for example, to minimize

$$\sum_{j \in \mathcal{J}_{k_0}} \int_{\mathcal{I}} (\hat{\phi}_j - \phi_{p(j)})^2,$$

then $\hat{p}$ is uniformly distributed on the set of all permutations of $\mathcal{J}_{k_0}$. From this, it may be shown, using properties (P) and letting $n \to \infty$ through values in $\mathcal{N}$, that (3.5) fails.

**4. Numerical properties.** This section summarizes the results of a Monte Carlo investigation of the finite-sample performance of the estimators $\hat{b}$ and $\tilde{b}$ discussed in Section 2. Samples of sizes $n = 100$ and 500 were generated from the model (2.1), with $\mathcal{I} = [0,1]$, $a = 0$ and the errors $\varepsilon_i$ distributed as normal $N(0, \sigma_\varepsilon^2)$, where $\sigma_\varepsilon = 0.5$ or 1. We took $b = \sum_{1 \leq j \leq 50} b_j \phi_j$ and $X = \sum_{1 \leq j \leq 50} \gamma_j Z_j \phi_j$, where (a) $b_1 = 0.3$ and $b_j = 4(-1)^{j+1} j^{-2}$ for $j > 1$, (b) the $\gamma_j$'s were deterministic coefficients, (c) $\phi_1 \equiv 1$ and $\phi_{j+1} = 2^{1/2} \cos(j\pi t)$ for $j \geq 1$ and (d) the $Z_j$'s were uniformly distributed on $[-3^{1/2}, 3^{1/2}]$. In particular, each $Z_j$ had zero mean and unit variance.

Two sets of the $\gamma_j$'s were used. In the first, $\gamma_j = (-1)^{j+1} j^{-\alpha/2}$, with $\alpha = 1.1, 1.5, 2$ or $4$. For these coefficients, the eigenvalues of the operator $K$ were $\kappa_j = j^{-\alpha}$ and were distinct. In the remainder of this section, we label these eigenvalues "well-spaced." In the second set, $\gamma_1 = 1$, $\gamma_j = 0.2(-1)^{j+1}(1 - 0.0001j)$ if $2 \leq j \leq 4$, and $\gamma_{5j+k} = 0.2(-1)^{5j+k+1} \{(5j)^{-\alpha/2} - 0.0001k\}$ for $j \geq 1$ and $0 \leq k \leq 4$. This set of $\gamma_j$'s generated blocks of $\kappa_j$'s that were nearly equal when $j$ was not too large and we refer to it as the "closely spaced" case. The theoretical arguments presented in Section 3 suggest that the performance of $\hat{b}$ can be poor in this setting.

All our results represent averages over 1000 Monte Carlo replications for each parameter setting. The quantities denoted by Bias$^2$, Var and MISE in Tables 1 and 2 are Monte Carlo approximations to integrated squared bias, integrated variance and mean integrated squared error, respectively, computed on a grid of 50 equally spaced points on $\mathcal{I}$. The values of $m$ and $\rho$, for given $n$, $\sigma_\varepsilon$, $\alpha$ and a given set of $\gamma_j$'s, were chosen to minimize MISE.

Table 1 shows that in the case of well-spaced eigenvalues, the MISE of $\hat{b}$ is smaller than that of $\tilde{b}$ for almost all values of the other design parameters. However, it follows from Table 2 that in the closely spaced case, the MISE of $\tilde{b}$ is nearly always smaller than that of $\hat{b}$. Thus, in terms of MISE, neither estimator dominates the other.

Both tables reveal that there is a general tendency for MISE to decrease as $\alpha$ increases. This does not contradict (3.5) or (3.10) since those results describe the behavior of MISE as a function of $n$ for fixed $\alpha$ and $\beta$, not the behavior of MISE as a function of $\alpha$ or $\beta$ for fixed $n$.

## 5. Derivation of Theorem 1.

5.1. *Proof of* (3.5). We begin by defining notation to be used in the proof. Given a sequence $c_n$ of positive constants, we shall use $O_p(c_n)$ and $o_p(c_n)$ to denote random variables $R_n$ and $r_n$, respectively, which satisfy

$$\lim_{D\to\infty} \limsup_{n\to\infty} \sup_{F\in\mathcal{F}} P_F(|R_n| > Dc_n) = 0,$$

$$\lim_{n\to\infty} \sup_{F\in\mathcal{F}} P_F(|r_n| > Dc_n) = 0 \qquad \text{for each } D > 0.$$

Similarly, a deterministic quantity $A_n = A_n(F)$, written as $A_n = O(c_n)$, will be understood to satisfy

$$\sup_{n\geq 1} c_n^{-1} \sup_{F\in\mathcal{F}} |A_n(F)| < \infty.$$

Next, we state subsidiary results concerning distances between the spectra of two operators. Let $L$ denote a general positive semidefinite linear operator as well as the kernel of that operator. Let the spectral decomposition of $L$ be

$$(5.1) \qquad L(u,v) = \sum_{j=1}^{\infty} \lambda_j \psi_j(u) \psi_j(v).$$

TABLE 1
*Results of Monte Carlo experiments for well-spaced eigenvalues*

| $\sigma_\varepsilon$ | $n$ | $\alpha$ | $\rho$ | $a_n$ | $\text{Bias}^2(\hat{b})$ | $\text{Bias}^2(\tilde{b})$ | $\text{Var}(\hat{b})$ | $\text{Var}(\tilde{b})$ | $\text{MISE}(\hat{b})$ | $\text{MISE}(\tilde{b})$ |
|---|---|---|---|---|---|---|---|---|---|---|
| 0.5 | 100 | 1.1 | 2 | 0.4 | 0.158 | 1.150 | 0.843 | 1.340 | 1.001 | 2.490 |
| | | 1.5 | 2 | 0.38 | 0.148 | 1.289 | 0.718 | 0.759 | 0.866 | 2.048 |
| | | 2.0 | 2 | 0.28 | 0.140 | 1.202 | 0.676 | 0.622 | 0.816 | 1.824 |
| | | 4.0 | 2 | 0.10 | 0.134 | 1.344 | 2.225 | 0.611 | 2.359 | 1.955 |
| | 500 | 1.1 | 3 | 0.28 | 0.016 | 0.717 | 0.236 | 0.480 | 0.251 | 1.197 |
| | | 1.5 | 3 | 0.22 | 0.015 | 0.663 | 0.254 | 0.364 | 0.269 | 1.027 |
| | | 2.0 | 2 | 0.12 | 0.139 | 0.416 | 0.146 | 0.441 | 0.285 | 0.857 |
| | | 4.0 | 2 | 0.032 | 0.139 | 0.460 | 0.409 | 0.493 | 0.548 | 0.953 |
| 1.0 | 100 | 1.1 | 2 | 1.0 | 0.161 | 2.709 | 2.034 | 1.203 | 2.195 | 3.913 |
| | | 1.5 | 2 | 0.75 | 0.149 | 2.401 | 2.221 | 1.019 | 2.370 | 3.420 |
| | | 2.0 | 2 | 0.50 | 0.139 | 2.047 | 2.395 | 1.034 | 2.534 | 3.081 |
| | | 4.0 | 1 | 0.25 | 3.257 | 2.302 | 0.501 | 0.788 | 3.758 | 3.090 |
| | 500 | 1.1 | 2 | 0.50 | 0.142 | 1.438 | 0.408 | 0.758 | 0.549 | 2.197 |
| | | 1.5 | 2 | 0.35 | 0.138 | 1.164 | 0.425 | 0.702 | 0.563 | 1.866 |
| | | 2.0 | 2 | 0.10 | 0.139 | 0.314 | 0.514 | 2.279 | 0.654 | 2.593 |
| | | 4.0 | 2 | 0.10 | 0.139 | 1.386 | 1.647 | 0.472 | 1.786 | 1.858 |

TABLE 2
*Results of Monte Carlo experiments for closely spaced eigenvalues*

| $\sigma_\varepsilon$ | $n$ | $\alpha$ | $m$ | $\rho$ | $\text{Bias}^2(\hat{b})$ | $\text{Bias}^2(\tilde{b})$ | $\text{Var}(\hat{b})$ | $\text{Var}(\tilde{b})$ | $\text{MISE}(\hat{b})$ | $\text{MISE}(\tilde{b})$ |
|---|---|---|---|---|---|---|---|---|---|---|
| 0.5 | 100 | 1.1 | 1 | 0.22 | 3.526 | 2.502 | 0.141 | 0.585 | 3.398 | 3.087 |
|  |  | 1.5 | 1 | 0.22 | 3.257 | 2.487 | 0.131 | 0.455 | 3.389 | 2.942 |
|  |  | 2.0 | 1 | 0.20 | 3.259 | 2.403 | 0.126 | 0.454 | 3.385 | 2.857 |
|  |  | 4.0 | 1 | 0.20 | 3.260 | 2.402 | 0.130 | 0.433 | 3.390 | 2.835 |
|  | 500 | 1.1 | 5 | 0.08 | 0.002 | 1.463 | 2.510 | 0.574 | 2.512 | 2.037 |
|  |  | 1.5 | 5 | 0.06 | 0.002 | 1.212 | 2.604 | 0.623 | 2.606 | 1.835 |
|  |  | 2.0 | 5 | 0.04 | 0.006 | 0.846 | 2.528 | 0.783 | 2.535 | 1.629 |
|  |  | 4.0 | 5 | 0.04 | 0.006 | 0.780 | 2.500 | 0.640 | 2.506 | 1.420 |
| 1.0 | 100 | 1.1 | 1 | 0.42 | 3.260 | 3.127 | 0.533 | 0.856 | 3.793 | 3.983 |
|  |  | 1.5 | 1 | 0.42 | 3.271 | 3.031 | 0.512 | 0.706 | 3.783 | 3.736 |
|  |  | 2.0 | 1 | 0.32 | 3.260 | 2.822 | 0.540 | 0.937 | 3.799 | 3.759 |
|  |  | 4.0 | 1 | 0.36 | 3.262 | 2.954 | 0.496 | 0.760 | 3.758 | 3.713 |
|  | 500 | 1.1 | 1 | 0.20 | 3.258 | 2.379 | 0.109 | 0.532 | 3.367 | 2.911 |
|  |  | 1.5 | 1 | 0.14 | 3.262 | 2.078 | 0.109 | 0.729 | 3.372 | 2.807 |
|  |  | 2.0 | 1 | 0.12 | 3.262 | 1.922 | 0.103 | 0.762 | 3.366 | 2.684 |
|  |  | 4.0 | 1 | 0.12 | 3.256 | 1.818 | 0.107 | 0.695 | 3.363 | 2.514 |

We assume that the terms are ordered in such a way that $\lambda_1 \geq \lambda_2 \geq \cdots \geq 0$. Given univariate functions $p$, $q$ and a symmetric bivariate function $M$, let $\|M\| = (\iint_{\mathcal{I}^2} M^2)^{1/2}$. Write $\int pq$ and $\int Mpq$ for

$$\int_\mathcal{I} p(u) q(u)\, du \quad \text{and} \quad \iint_{\mathcal{I}^2} M(u,v) p(u) p(v)\, du\, dv,$$

respectively. Further, denote by $\int Mp$ the function of which the value at $u$ is $\int_\mathcal{I} M(u,v) p(v)\, dv$ and define $\delta_j = \min_{1 \leq k \leq j} (\kappa_k - \kappa_{k+1})$.

The following pair of results may be derived from theory developed by Bhatia, Davis and McIntosh [1]:

(5.2) $$\sup_{j \geq 1} |\kappa_j - \lambda_j| \leq \|K - L\|, \qquad \sup_{j \geq 1} \delta_j \|\phi_j - \psi_j\| \leq 8^{1/2} \|K - L\|.$$

In framing the second bound here, we use the convention that $\int \psi_j \phi_j \geq 0$. This determines the sign of $\psi_j$ in those cases where choice of sign has an impact on the validity of (5.2).

The following lemma will be proven in Section 5.2:

LEMMA 5.1. *If we are able to write $\psi_j - \phi_j = \chi_j + \Delta_j$ for functions $\chi_j$ and $\Delta_j$, then*

(5.3)
$$\left| (\kappa_j - \lambda_j) \left(1 + \int \chi_j \phi_j\right) - \int (K - L)(\phi_j + \chi_j) \phi_j \right|$$
$$\leq \|\Delta_j\| \left\{ |\kappa_j - \lambda_j| + \left\| \int (K - L) \phi_j \right\| \right\}.$$

*Furthermore, if* $\inf_{k:\,k\neq j}|\lambda_j - \kappa_k| > 0$, *then*

$$(5.4) \quad \psi_j - \phi_j = \sum_{k:\,k\neq j}(\lambda_j - \kappa_k)^{-1}\phi_k \int (L-K)\psi_j\phi_k + \phi_j \int (\psi_j - \phi_j)\phi_j.$$

Put $\widehat{\Delta} = \|\widehat{K} - K\|$ and define the event $\mathcal{E}_m$ by

$$\mathcal{E}_m = \mathcal{E}_m(n) = \{\tfrac{1}{2}\kappa_m \geq \widehat{\Delta}\}.$$

That is, $\mathcal{E}_m$ denotes the set of all realizations such that for sample size $n$, $\tfrac{1}{2}\kappa_m \geq \widehat{\Delta}$. Below, when we say that a bound is valid when $\mathcal{E}_m$ holds, this should be interpreted as stating that the bound is valid for all realizations for which $\tfrac{1}{2}\kappa_m \geq \widehat{\Delta}$. It is not a statement that relates to a conditioning argument in the sense that conditioning is usually interpreted in probability theory.

Write $\hat{b}_j = \check{b}_j + \hat{\kappa}_j^{-1}(S_{j1} + S_{j2} + S_{j3})$, where $\hat{\kappa}_j \check{b}_j = \int g\phi_j$, $S_{j1} = \int (\hat{g} - g)\phi_j$, $S_{j2} = \int g(\hat{\phi}_j - \phi_j)$ and $S_{j3} = \int (\hat{g} - g)(\hat{\phi}_j - \phi_j)$. In this notation,

$$(5.5) \quad \begin{aligned} \sum_{j=1}^m (\hat{b}_j - \check{b}_j)^2 &\leq 3 \sum_{j=1}^m \hat{\kappa}_j^{-2}(S_{j1}^2 + S_{j2}^2 + S_{j3}^2) \\ &\leq 12 \sum_{j=1}^m \kappa_j^{-2}(S_{j1}^2 + S_{j2}^2 + S_{j3}^2) \\ &\leq 12 \sum_{j=1}^m \kappa_j^{-2}(S_{j1}^2 + S_{j2}^2) + 12\|\hat{g} - g\|^2 \sum_{j=1}^m \kappa_j^{-2}\|\hat{\phi}_j - \phi_j\|^2, \end{aligned}$$

where the first inequality holds universally; the second inequality, obtained using the first part of (5.2), is valid provided the event $\mathcal{E}_m$ holds; and the third inequality employs the bound $|S_{j3}| \leq \|\hat{g} - g\|\,\|\hat{\phi}_j - \phi_j\|$.

Note that provided $\mathcal{E}_m$ holds, we have

$$(5.6) \quad \begin{aligned} \sum_{j=1}^m (\check{b}_j - b_j)^2 &= \sum_{j=1}^m \left(\frac{\hat{\kappa}_j - \kappa_j}{\hat{\kappa}_j \kappa_j}\right)^2 \left(\int g\phi_j\right)^2 \leq 4 \sum_{j=1}^m \left(\frac{\hat{\kappa}_j - \kappa_j}{\kappa_j^2}\right)^2 \left(\int g\phi_j\right)^2 \\ &= 4 \sum_{j=1}^m \kappa_j^{-2} b_j^2 (\hat{\kappa}_j - \kappa_j)^2. \end{aligned}$$

Define $\widehat{\Delta}_j = \|\int (\widehat{K} - K)\phi_j\|$. Using (5.3) with $\chi_j \equiv 0$ and then applying both parts of (5.2), we obtain

$$(5.7) \quad \begin{aligned} \left|\hat{\kappa}_j - \kappa_j - \int (\widehat{K} - K)\phi_j\phi_j\right| &\leq \|\hat{\phi}_j - \phi_j\|(|\hat{\kappa}_j - \kappa_j| + \widehat{\Delta}_j) \\ &\leq \delta_j^{-1}\widehat{\Delta}(\widehat{\Delta} + \widehat{\Delta}_j). \end{aligned}$$

Combining (5.6) and (5.7) and defining $\widehat{\Delta}_{jj} = |\int (\widehat{K} - K)\phi_j\phi_j|$, we deduce that if $\mathcal{E}_m$ holds, then

$$(5.8) \quad \sum_{j=1}^{m}(\check{b}_j - b_j)^2 \leq 8\sum_{j=1}^{m}\kappa_j^{-2}b_j^2\widehat{\Delta}_{jj}^2 + 16\widehat{\Delta}^2\sum_{j=1}^{m}(\delta_j\kappa_j)^{-2}b_j^2(\widehat{\Delta}^2 + \widehat{\Delta}_j^2).$$

We shall prove in Section 5.3 that under the conditions of the theorem,

$$(5.9) \quad E(\widehat{\Delta}^2) + E(\widehat{\Delta}_j^2) = O(n^{-1}), \qquad E(\widehat{\Delta}_{jj}^2) = O(n^{-1}\kappa_j^2),$$

uniformly in $j$. In particular, (5.9) entails $\widehat{\Delta} = O_p(n^{-1/2})$. Now, (3.2) and (3.4) imply that $n^{1/2}\kappa_m \to \infty$ as $n \to \infty$, so the first part of (5.9) implies that $P(\mathcal{E}_m) \to 1$. Therefore, since the result (3.5) that we wish to prove relates only to probabilities of differences (not to moments of differences), it suffices to work with bounds that are established under the assumption that $\mathcal{E}_m$ holds, since the contrary case contributes only $o(1)$ to the probability on the left-hand side of (3.5).

In our arguments below, we shall use the property $\widehat{\Delta} = O_p(n^{-1/2})$ without further reference. Now, the conditions in Theorem 1 imply that $\delta_j^{-1} \leq C_1 j^{\alpha+1}$, whence it follows that

$$(5.10) \quad \begin{aligned} n^{-1}\sum_{j=1}^{m}\kappa_j^{-2}b_j^2\kappa_j^2 &\leq C_2 n^{-1}\sum_{j=1}^{m} j^{-2\beta} \leq C_3 n^{-1}, \\ n^{-2}\sum_{j=1}^{m}(\delta_j\kappa_j)^{-2}b_j^2 &\leq C_4 n^{-2}\sum_{j=1}^{m} j^{4\alpha-2\beta+2} \leq C_5 n^{-2}s(n), \end{aligned}$$

where $C_1, \ldots, C_5$ are positive constants and $s(n)$ equals $n^{(4\alpha-2\beta+3)/(\alpha+2\beta)}$ if the exponent is strictly positive, equals $1 + \log n$ if the exponent vanishes and equals 1 otherwise. Combining (5.8)–(5.10), we deduce that

$$(5.11) \quad \sum_{j=1}^{m}(\check{b}_j - b_j)^2 = O_p\{n^{-1} + n^{-2}s(n)\} = o_p(n^{-(2\beta-1)/(\alpha+2\beta)}).$$

Observe, too, that

$$(5.12) \quad \begin{aligned} \int_{\mathcal{I}}\left\{\sum_{j=1}^{m}b_j\hat{\phi}_j(u) - b(u)\right\}^2 du &\leq 2\int_{\mathcal{I}}\left[\sum_{j=1}^{m}b_j\{\hat{\phi}_j(u) - \phi_j(u)\}\right]^2 du \\ &\quad + 2\sum_{j=m+1}^{\infty}b_j^2 \\ &\leq 2m\sum_{j=1}^{m}b_j^2\|\hat{\phi}_j - \phi_j\|^2 + 2\sum_{j=m+1}^{\infty}b_j^2. \end{aligned}$$

Combining (5.5), (5.11) and (5.12), we find that

$$\int (\hat{b} - b)^2 \leq 3 \sum_{j=1}^{m} (\hat{b}_j - \check{b}_j)^2 + 3 \sum_{j=1}^{m} (\check{b}_j - b_j)^2$$

$$+ 3 \sum_{j=1}^{m} \int \left( \sum_{j=1}^{m} b_j \hat{\phi}_j - b \right)^2$$

(5.13)
$$\leq 36 \sum_{j=1}^{m} \kappa_j^{-2} (S_{j1}^2 + S_{j2}^2)$$

$$+ 36 \sum_{j=1}^{m} (m b_j^2 + \|\hat{g} - g\|^2 \kappa_j^{-2}) \|\hat{\phi}_j - \phi_j\|^2$$

$$+ 6 \sum_{j=m+1}^{\infty} b_j^2 + o_p(n^{-(2\beta-1)/(\alpha+2\beta)}).$$

Simple moment calculations show that $E\|\hat{g} - g\|^2 = O(n^{-1})$ and, clearly, $\sum_{j \geq m+1} b_j^2 = O(n^{-(2\beta-1)/(\alpha+2\beta)})$. It will be proved in Section 5.3 that

(5.14) $$E(S_{j1}^2) = O(n^{-1} \kappa_j),$$

whence it follows that $\sum_{j \leq m} \kappa_j^{-2} S_{j1}^2 = O_p(n^{-(2\beta-1)/(\alpha+2\beta)})$. Combining these results and (5.13), we see that (3.5) will follow if we prove that

(5.15) $$\sum_{j=1}^{m} j^{2\alpha} \left\{ \int g(\hat{\phi}_j - \phi_j) \right\}^2 + \sum_{j=1}^{m} (m j^{-2\beta} + n^{-1} j^{2\alpha}) \|\hat{\phi}_j - \phi_j\|^2$$
$$= O_p(n^{-(2\beta-1)/(\alpha+2\beta)}).$$

Derivation of this property requires bounds on $\hat{\phi}_j - \phi_j$, which we now discuss.

Take $L = \widehat{K}$, $\lambda_j = \hat{\kappa}_j$ and $\psi_j = \hat{\phi}_j$ in Lemma 5.1. Formula (5.4) yields $\|\hat{\phi}_j - \phi_j\|^2 = \hat{u}_j^2 + \hat{v}_j^2$, where

$$\hat{u}_j^2 = \sum_{k: k \neq j} (\hat{\kappa}_j - \kappa_k)^{-2} \left\{ \int (\widehat{K} - K) \hat{\phi}_j \phi_k \right\}^2$$

and $\hat{v}_j^2 = \{\int (\hat{\phi}_j - \phi_j) \phi_j\}^2$. Now, $\hat{u}_j$ equals the length of the projection of $\hat{\phi}_j - \phi_j$ into the plane perpendicular to $\phi_j$; hence, it also equals the projection of $\hat{\phi}_j$ into that plane. Also, $\int \hat{\phi}_j \phi_j$ equals the length of the projection of $\hat{\phi}_j$ onto $\phi_j$. Therefore, by Pythagoras' Theorem, $(\int \hat{\phi}_j \phi_j)^2 + \hat{u}_j^2 = \|\hat{\phi}_j\|^2 = 1$, whence it follows that $\int \hat{\phi}_j \phi_j = (1 - \hat{u}_j^2)^{1/2}$. Hence,

$$\hat{v}_j^2 = \left( 1 - \int \hat{\phi}_j \phi_j \right)^2 = \{1 - (1 - \hat{u}_j^2)^{1/2}\}^2 = 2\{1 - (1 - \hat{u}_j^2)^{1/2}\} - \hat{u}_j^2,$$

which implies that

(5.16) $$\|\hat{\phi}_j - \phi_j\|^2 = 2\{1 - (1 - \hat{u}_j^2)^{1/2}\} \leq 2\hat{u}_j^2.$$

Let $C > 0$ and define

$$\mathcal{F}_m = \mathcal{F}_m(n) = \{(\hat{\kappa}_j - \kappa_k)^{-2} \leq 2(\kappa_j - \kappa_k)^{-2} \leq Cn^{2(\alpha+1)/(\alpha+2\beta)}\},$$

that is, the set of realisations such that, for sample size $n$, $(\hat{\kappa}_j - \kappa_k)^{-2} \leq 2(\kappa_j - \kappa_k)^{-2} \leq Cn^{2(\alpha+1)/(\alpha+2\beta)}$. Observe that

(5.17) $$\left\{\int(\widehat{K} - K)\hat{\phi}_j\phi_k\right\}^2 \leq 2\left\{\int(\widehat{K} - K)\phi_j\phi_k\right\}^2 + 2\hat{w}_{jk}^2,$$

where $\hat{w}_{jk}^2 = \{\int(\widehat{K} - K)(\hat{\phi}_j - \phi_j)\phi_k\}^2$. Note, too, that uniformly in $1 \leq j \leq m$,

$$\max(\kappa_j - \kappa_{j+1}, \kappa_{j-1} - \kappa_j) \geq C_1 j^{-(\alpha+1)} \geq C_2 n^{-(\alpha+1)/(\alpha+2\beta)},$$

where $C_1, C_2$ denote positive constants, and that since $\beta > \frac{1}{2}\alpha + 1$, it follows that $n^{-1/2} = o(n^{-(\alpha+1)/(\alpha+2\beta)})$. These properties, and the fact that $|\hat{\kappa}_j - \kappa_j| \leq \widehat{\Delta} = O_p(n^{-1/2})$, imply that if the constant $C$ in the definition of $\mathcal{F}_m$ is chosen to be sufficiently large, then $P(\mathcal{F}_m) \to 1$ as $n \to \infty$. Therefore, as in the case of $\mathcal{E}_n$, since (3.5) relates only to probabilities of differences, it suffices to work with bounds that are established under the assumption that $\mathcal{F}_m$ holds. In this case,

(5.18) $$\sum_{k:\,k\neq j}(\hat{\kappa}_j - \kappa_k)^{-2}\hat{w}_{jk}^2 \leq Cn^{2(\alpha+1)/(\alpha+2\beta)}\sum_{k=1}^{\infty}\hat{w}_{jk}^2.$$

Using Parseval's identity and the Cauchy–Schwarz inequality, we may prove that

(5.19) $$\sum_{k=1}^{\infty}\hat{w}_{jk}^2 = \int_{\mathcal{I}}\left[\int_{\mathcal{I}}(\widehat{K} - K)(u,v)(\hat{\phi}_j - \phi_j)(v)\,dv\right]^2 du \leq \widehat{\Delta}^2\|\hat{\phi}_j - \phi_j\|^2.$$

Combining (5.17)–(5.19), we deduce that provided $\mathcal{F}_m$ holds, we have

$$\hat{u}_j^2 \leq 2\sum_{k:\,k\neq j}(\hat{\kappa}_j - \kappa_k)^{-2}\left\{\int(\widehat{K} - K)\phi_j\phi_k\right\}^2 + 2Cn^{2(\alpha+1)/(\alpha+2\beta)}\widehat{\Delta}^2\|\hat{\phi}_j - \phi_j\|^2.$$

Substituting into (5.16), we find that

(5.20) $$\begin{aligned}(1 - 4Cn^{2(\alpha+1)/(\alpha+2\beta)}\widehat{\Delta}^2)\|\hat{\phi}_j - \phi_j\|^2 \\ \leq 4\sum_{k:\,k\neq j}(\hat{\kappa}_j - \kappa_k)^{-2}\left\{\int(\widehat{K} - K)\phi_j\phi_k\right\}^2.\end{aligned}$$

Recall that $\widehat{\Delta} = O_p(n^{-1/2})$ and observe that since $\beta > \frac{1}{2}\alpha + 1$, we have $n^{2(\alpha+1)/(\alpha+2\beta)} \cdot n^{-1} \to 0$. Therefore, noting that $P(\mathcal{F}_m) \to 1$, we deduce that (5.20) implies

$$
\begin{aligned}
(5.21) \quad \|\hat{\phi}_j - \phi_j\|^2 &\leq 4\{1 + o_p(1)\} \sum_{k:k\neq j} (\hat{\kappa}_j - \kappa_k)^{-2} \left\{ \int (\widehat{K} - K)\phi_j\phi_k \right\}^2 \\
&\leq 8\{1 + o_p(1)\} \sum_{k:k\neq j} (\kappa_j - \kappa_k)^{-2} \left\{ \int (\widehat{K} - K)\phi_j\phi_k \right\}^2,
\end{aligned}
$$

where the $o_p(1)$ terms are of that order uniformly in $1 \leq j \leq m$. We shall show in Section 5.3 that

$$
(5.22) \quad n \sum_{k:k\neq j} (\kappa_j - \kappa_k)^{-2} E\left\{ \int (\widehat{K} - K)\phi_j\phi_k \right\}^2 = O(j^2),
$$

uniformly in $1 \leq j \leq m$. Results (5.21) and (5.22) together imply that

$$
(5.23) \quad \sum_{j=1}^m (mj^{-2\beta} + n^{-1}j^{2\alpha})\|\hat{\phi}_j - \phi_j\|^2 = O_p(mn^{-1} + m^{2\alpha+3}n^{-2})
$$
$$
= o_p(n^{-(2\beta-1)/(\alpha+2\beta)}).
$$

Next, observe that

$$
\begin{aligned}
(5.24) \quad \int g(\hat{\phi}_j - \phi_j) &= \sum_{k:k\neq j} g_k(\hat{\kappa}_j - \kappa_k)^{-1} \int (\widehat{K} - K)\hat{\phi}_j\phi_k \\
&\quad + g_j \int (\hat{\phi}_j - \phi_j)\phi_j \\
&= T_{j1} + T_{j2} + T_{j3} + T_{j4},
\end{aligned}
$$

where

$$
T_{j1} = \sum_{k:k\neq j} g_k (\kappa_j - \kappa_k)^{-1} \int (\widehat{K} - K)\phi_j\phi_k,
$$

$$
T_{j2} = \sum_{k:k\neq j} g_k\{(\hat{\kappa}_j - \kappa_k)^{-1} - (\kappa_j - \kappa_k)^{-1}\} \int (\widehat{K} - K)\phi_j\phi_k,
$$

$$
T_{j3} = \sum_{k:k\neq j} g_k(\hat{\kappa}_j - \kappa_k)^{-1} \int (\widehat{K} - K)(\hat{\phi}_j - \phi_j)\phi_k
$$

and $T_{j4} = g_j \int (\hat{\phi}_j - \phi_j)\phi_j$. Let $C_1, C_2, \ldots$ denote positive constants. Since $|g_k| \leq C_1 k^{-(\alpha+\beta)}$, then if $\mathcal{F}_m$ holds, we have

$$
T_{j2}^2 \leq C_2 \left\{ \sum_{k:k\neq j} k^{-(\alpha+\beta)} \frac{|\hat{\kappa}_j - \kappa_j|}{(\kappa_j - \kappa_k)^2} \left| \int (\widehat{K} - K)\phi_j\phi_k \right| \right\}^2
$$

$$\leq C_3 \bigg\{ \sum_{k:\,k\neq j} k^{-2(\alpha+\beta)} \frac{(\hat{\kappa}_j - \kappa_j)^2}{(\kappa_j - \kappa_k)^4} \bigg\} \bigg[ \sum_{k=1}^{\infty} \bigg\{ \int (\widehat{K} - K)\phi_j \phi_k \bigg\}^2 \bigg]$$

$$= C_4 (\hat{\kappa}_j - \kappa_j)^2 \widehat{\Delta}_j^2 \sum_{k:\,k\neq j} k^{-2(\alpha+\beta)} (\kappa_j - \kappa_k)^{-4}.$$

Now,

$$\sum_{k=2j}^{\infty} k^{-2(\alpha+\beta)} (\kappa_j - \kappa_k)^{-4} \leq C_5 \kappa_j^{-4} \sum_{k=2j}^{\infty} k^{-2(\alpha+\beta)} \leq C_6 j^{2\alpha - 2\beta + 1},$$

$$\sum_{k=1}^{j/2} k^{-2(\alpha+\beta)} (\kappa_j - \kappa_k)^{-4} \leq C_7 \sum_{k=1}^{j/2} k^{-2(\alpha+\beta)} \kappa_k^{-4} \leq C_8 \sum_{k=1}^{j/2} k^{2\alpha - 2\beta}$$

$$\leq C_9 \begin{cases} 1, & \text{if } \alpha + \tfrac{1}{2} < \beta, \\ 1 + \log j, & \text{if } \alpha + \tfrac{1}{2} = \beta, \\ j^{2\alpha - 2\beta + 1}, & \text{if } \alpha + \tfrac{1}{2} > \beta, \end{cases}$$

$$\sum_{k=j/2}^{2j} k^{-2(\alpha+\beta)} (\kappa_j - \kappa_k)^{-4} \leq \sum_{k=j/2}^{2j} k^{-2(\alpha+\beta)} (j/\kappa_j)^4 (1 + |j-k|)^{-4}$$

$$\leq C_{10} j^{2\alpha - 2\beta + 4}.$$

Therefore,

$$(5.25) \qquad \sum_{k:\,k\neq j} k^{-2(\alpha+\beta)} (\kappa_j - \kappa_k)^{-4} \leq C_{11} (1 + j^{2\alpha - 2\beta + 4} + \log j),$$

whence, using (5.2) and (5.9), we have

$$\sum_{j=1}^{m} j^{2\alpha} T_{j2}^2 \leq C_{12} \sum_{j=1}^{m} (\hat{\kappa}_j - \kappa_j)^2 \widehat{\Delta}_j^2 (j^{2\alpha} \log n + j^{4\alpha - 2\beta + 4})$$

$$(5.26) \qquad = O_p \bigg\{ n^{-1} \sum_{j=1}^{m} E(\widehat{\Delta}_j^2)(j^{2\alpha} \log n + j^{4\alpha - 2\beta + 4}) \bigg\}$$

$$= O_p \bigg\{ n^{-2} \sum_{j=1}^{m} (j^{2\alpha} \log n + j^{4\alpha - 2\beta + 4}) \bigg\}$$

$$= O_p \{ n^{-2} (m^{2\alpha + 1} \log n + m^{4\alpha - 2\beta + 5}) \} = o_p(n^{-(2\beta - 1)/(\alpha + 2\beta)}).$$

If $\mathcal{F}_m$ holds, then

$$
\begin{aligned}
|T_{j3}| &\leq C_{13} \sum_{k:\,k\neq j} k^{-(\alpha+\beta)}|\kappa_j - \kappa_k|^{-1}\left|\int(\widehat{K} - K)(\hat{\phi}_j - \phi_j)\phi_k\right| \\
&\leq C_{14} \sum_{k:\,k\neq j} k^{-(\alpha+\beta)}|\kappa_j - \kappa_k|^{-1}\|\hat{\phi}_j - \phi_j\| \\
&\quad \times \int |\phi_k(u)|\left[\int\{\widehat{K}(u,v) - K(u,v)\}^2\,dv\right]^{1/2} du \\
&\leq C_{15}\widehat{\Delta}\|\hat{\phi}_j - \phi_j\|\sum_{k:\,k\neq j} k^{-(\alpha+\beta)}|\kappa_j - \kappa_k|^{-1} \\
&\leq C_{16}\widehat{\Delta}\|\hat{\phi}_j - \phi_j\|,
\end{aligned}
\quad (5.27)
$$

where the last inequality follows using the argument leading to (5.25). From (5.27), using (5.21) and (5.22), it may be shown that

$$
(5.28)\quad \sum_{j=1}^{m} j^{2\alpha} T_{j3}^2 = O_p\left(n^{-2}\sum_{j=1}^{m} j^{2\alpha+2}\right) = o_p(n^{-(2\beta-1)/(\alpha+2\beta)}).
$$

More simply,

$$
(5.29)\quad \sum_{j=1}^{m} j^{2\alpha} T_{j4}^2 \leq C_{17}\sum_{j=1}^{m} j^{-2\beta}\|\hat{\phi}_j - \phi_j\|^2 = O_p(n^{-1}).
$$

Combining (5.24), (5.26), (5.28) and (5.29), we deduce that

$$
(5.30)\quad \sum_{j=1}^{m} j^{2\alpha}\left\{\int g(\hat{\phi}_j - \phi_j)\right\}^2 \leq 4\sum_{j=1}^{m} j^{2\alpha} T_{j1}^2 + o_p(n^{-(2\beta-1)/(\alpha+2\beta)}).
$$

We shall prove in Section 5.3 that

$$
(5.31)\quad \sum_{j=1}^{m} j^{2\alpha} E(T_{j1}^2) = O(n^{-(2\beta-1)/(\alpha+2\beta)}).
$$

The desired result (5.15) follows from (5.23), (5.30) and (5.31). This completes the proof of (3.5).

5.2. *Proof of Lemma* 5.1. To derive (5.3), observe that on subtracting the expansions of $K$ and $L$ in (2.1) and (5.1), respectively, we obtain an expansion of $K - L$. Multiplying both sides of this by $\psi_j(u)\phi_j(v)$ and integrating over $u$ and $v$, we deduce that

$$
(5.32)\quad (\kappa_j - \lambda_j)\int \psi_j\phi_j - \int(K-L)\psi_j\phi_j = 0.
$$

Since $\psi_j = \phi_j + \chi_j + \Delta_j$, we have

$$
(5.33)\quad \left|\int \psi_j\phi_j - 1 - \int \chi_j\phi_j\right| = \left|\int \Delta_j\phi_j\right| \leq \|\Delta_j\|,
$$

$$
\begin{aligned}
\left|\int (K-L)(\psi_j - \phi_j - \chi_j)\phi_j\right|^2 & \\
(5.34) \quad = \left|\int (K-L)\Delta_j \phi_j\right|^2 & \\
\leq \left(\int \Delta_j^2\right) \int_{\mathcal{I}} \left[\int_{\mathcal{I}} \{K(u,v) - L(u,v)\}\phi_j(u)\,du\right]^2 dv.
\end{aligned}
$$

Result (5.3) follows from (5.32)–(5.34).

The expansions of $K$ and $L$ in (2.1) and (5.1) may be used to prove that

$$\lambda_j(\psi_j - \phi_j) = \int K(\psi_j - \phi_j) + \int (L-K)\psi_j - (\lambda_j - \kappa_j)\phi_j.$$

Multiplying both sides by $\phi_k$ and integrating, we deduce that

$$\lambda_j \int (\psi_j - \phi_j)\phi_k = \kappa_k \int (\psi_j - \phi_j)\phi_k + \int (L-K)\psi_j\phi_k - (\lambda_j - \kappa_j)\delta_{jk},$$

where $\delta_{jk}$ denotes the Kronecker delta. Equivalently, provided $\lambda_j \neq \kappa_k$, we have

$$\int (\psi_j - \phi_j)\phi_k = (\lambda_j - \kappa_k)^{-1} \int (L-K)\psi_j\phi_k - \delta_{jk}.$$

Result (5.4) follows from this formula and the fact that

$$\psi_j - \phi_j = \sum_{k=1}^{\infty} \phi_k \int (\psi_j - \phi_j)\phi_k.$$

5.3. *Proofs of* (5.9), (5.14), (5.22) *and* (5.31). Direct calculation shows that $E(\widehat{K} - K)^2 = O(n^{-1})$, uniformly on $\mathcal{I} \times \mathcal{I}$. It follows that $E(\widehat{\Delta}^2) = O(n^{-1})$. Note, too, that by Parseval's identity, $\widehat{\Delta}^2 = \sum_j \widehat{\Delta}_j^2$ and so $\sup_j E(\widehat{\Delta}_j^2) = O(n^{-1})$.

This gives the first part of (5.9). To derive the second part, assume without loss of generality that $E(X) = 0$ and observe that

$$(5.35) \quad \int (\widehat{K} - K)\phi_j \phi_k = n^{-1} \sum_{i=1}^{n} (\xi_{ij}\xi_{ik} - E\xi_j\xi_k) - \bar{\xi}_j \bar{\xi}_k,$$

where $\xi_{ij} = \int X_i \phi_j$, $\bar{\xi}_j = n^{-1} \sum_i \xi_{ij}$ and $\xi_j$ denotes a generic $\xi_{ij}$. Therefore, using the fact that $E(\xi_j^4) \leq C_1 (E\xi_j^2)^2$, where $C_1 > 0$ does not depend on $j$, we deduce that

$$E(\widehat{\Delta}_{jj}^2) = E\left\{n^{-1}\sum_{i=1}^{n}(\xi_{ij}^2 - E\xi_j^2) - \bar{\xi}_j^2\right\}^2 \leq n^{-1}C_2(E\xi_j^2)^2 = n^{-1}C_2\kappa_j^2,$$

where $C_2$ does not depend on $j$. This implies the second part of (5.9).

To prove (5.14), observe that

$$\int (\hat{g} - g)\phi_j = n^{-1} \sum_{i=1}^{n} \left\{ \xi_{ij} \int bX_i - \bar{\xi}_j \int b\bar{X} - E\left(\xi_{ij} \int bX_i\right) \right\}$$
$$+ n^{-1} \sum_{i=1}^{n} (\xi_{ij}\varepsilon_i - \bar{\xi}_j \bar{\varepsilon}),$$

where $\bar{\varepsilon} = n^{-1} \sum_i \varepsilon_i$. It may thus be proved that

$$nE\left\{ \int (\hat{g} - g)\phi_j \right\}^2 \leq C_3 \left\{ \text{var}\left(\xi_j \int bX\right) + \text{var}(\xi_j) \right\} \leq C_4 (E\xi_j^4)^{1/2} \leq C_5 \kappa_j,$$

which implies (5.14).

To obtain (5.22), note that by (5.35) and the fact that $E(\xi_j^4) \leq C_1(E\xi_j^2)^2$, we have

$$nE\left\{ \int (\widehat{K} - K)\phi_j \phi_k \right\}^2 \leq C_6 E(\xi_j^2 \xi_k^2) \leq C_7 (E\xi_j^4 \cdot E\xi_k^4)^{1/2} \leq C_8 \kappa_j \kappa_k,$$

uniformly in $j$ and $k$. Result (5.22) follows directly on substitution and employing the argument leading to (5.25).

Again using (5.35), we have

$$T_{j1} = n^{-1} \sum_{i=1}^{n} \sum_{k:k\neq j} g_k(\kappa_j - \kappa_k)^{-1} \{\xi_{ij}\xi_{ik} - E(\xi_{ij}\xi_{ik}) - \bar{\xi}_j \bar{\xi}_k\},$$

from which it may be proved that since $E(|\xi_{k_1} \ldots \xi_{k_4}|) \leq \prod_\ell (E\xi_{k_\ell}^4)^{1/4}$,

$$nE(T_{j1}^2) \leq C_9 E\left\{ \xi_j \sum_{k:k\neq j} g_k(\kappa_j - \kappa_k)^{-1} \xi_k \right\}^2$$
$$\leq C_{10}(E\xi_j^4)^{1/2} \left\{ \sum_{k_1:k_1\neq j} \cdots \sum_{k_4:k_4\neq j} E(|\xi_{k_1} \ldots \xi_{k_4}|) \right.$$
$$\left. \times \prod_{\ell=1}^{4} |g_{k_\ell}(\kappa_j - \kappa_{k_\ell})^{-1}| \right\}^{1/2}$$
$$\leq C_{11}\kappa_j \left( \sum_{k:k\neq j} \left|\frac{g_k \kappa_k}{\kappa_j - \kappa_k}\right| \right)^4 \leq C_{12}\kappa_j,$$

uniformly in $j$. Therefore, $\sum_{j\leq m} j^{2\alpha} E(T_{j1}^2) \leq C_{10} n^{-1} m^{\alpha+1}$, which implies (5.31).

5.4. *Proof of* (3.6). Let $\mathcal{I} \equiv [0,1]$, $\phi_1 \equiv 1$ and $\phi_{j+1}(t) = 2^{-1/2}\cos(j\pi t)$ for $j \geq 1$. Put $b_j = \theta_j j^{-\beta}$ for $L_{n+1} \leq j \leq 2L_n$ and $b_j = 0$ otherwise, where $L_n$ denotes the integer part of $n^{1/(2\beta+1)}$ and each $\theta_j$ is either 0 or 1. Let $\kappa_j = j^{-\alpha}$ and write $Z_1, Z_2, \ldots$ for independent random variables, all with the uniform distribution on $[-3^{1/2}, 3^{1/2}]$. Note that $E(Z_j) = 0$, $E(Z_j^2) = 1$ and that the $Z_j$'s are observable if $X$ is observable since $Z_j = j^{\alpha/2} \int_\mathcal{I} X \phi_j$.

Set $X = \sum_j j^{-\alpha/2} Z_j \phi_j$ and

$$Y = \int_\mathcal{I} bX + \varepsilon = \sum_{j=L_n+1}^{2L_n} \theta_j j^{-(\alpha+2\beta)/2} Z_j + \varepsilon,$$

where the error $\varepsilon$ is taken to be Gaussian with zero mean. Then we may write $b = \sum_{L_n+1 \leq j \leq 2L_n} \theta_j j^{-\beta} \phi_j$ and if $\bar{b}$ is an estimator of $b$, it follows that

$$(5.36) \qquad \bar{\theta}_j = j^\beta \int_\mathcal{I} \bar{b} \phi_j$$

is an estimator of $\theta_j$. An argument based on the Neyman–Pearson lemma shows that

$$\lim_{n \to \infty} \inf_{L_n+1 \leq j \leq 2L_n} \inf_{\bar{\theta}_j} \sup{}^* E(\bar{\theta}_j - \theta_j)^2 > 0,$$

where $\sup^*$ denotes the supremum over all $2^{L_n}$ different distributions of $(X, Y)$ obtained by taking different choices of $\theta_{L_n+1}, \ldots, \theta_{2L_n}$, and $\inf_{\bar{\theta}_j}$ represents the infimum over all measurable functions $\bar{\theta}_j$ of the data. Therefore, if an estimator $\check{b}$ is given and $\check{\theta}_{L_n+1}, \ldots, \check{\theta}_{2L_n}$ are the respective estimators of $\theta_{L_n+1}, \ldots, \theta_{2L_n}$ obtained by substituting $\check{b}$ for $\bar{b}$ in (5.36), then for constants $D_1, D_2 > 0$ which do not depend on the choice of the measurable function $\check{b}$,

$$\sup{}^* \int_\mathcal{I} E_F(\check{b} - b)^2 = \sup{}^* \sum_{j=L_n+1}^{2L_n} j^{-2\beta} E_F(\check{\theta}_j - \theta_j)^2$$

$$\geq D_1 \sum_{j=L_n+1}^{2L_n} j^{-2\beta} \geq D_2 n^{-(2\beta-1)/(\alpha+2\beta)}.$$

This proves (3.6).

## REFERENCES

[1] BHATIA, R., DAVIS, C. and MCINTOSH, A. (1983). Perturbation of spectral subspaces and solution of linear operator equations. *Linear Algebra Appl.* **52/53** 45–67. MR0709344
[2] BLUNDELL, R. and POWELL, J. L. (2003). Endogeneity in nonparametric and semiparametric regression models. In *Advances in Economics and Econometrics: Theory and Applications* (M. Dewatripont, L. P. Hansen and S. J. Turnovsky, eds.) **2** 312–357. Cambridge Univ. Press.
[3] CARDOT, H. and SARDA, P. (2003). Linear regression models for functional data. Available at www.quantlet.org/hizirjsp/sarda-cardot/sarda-cardot.pdf

[4] CARROLL, R. J., RUPPERT, D. and STEFANSKI, L. A. (1995). *Measurement Error in Nonlinear Models*. Chapman and Hall, London. MR1630517

[5] CAVALIER, L., GOLUBEV, G. K., PICARD, D. and TSYBAKOV, A. B. (2002). Oracle inequalities for inverse problems. *Ann. Statist.* **30** 843–874. MR1922543

[6] CHIOU, J.-M., MÜLLER, H.-G. and WANG, J.-L. (2003). Functional quasi-likelihood regression models with smooth random effects. *J. R. Stat. Soc. Ser. B Stat. Methodol.* **65** 405–423. MR1983755

[7] CUEVAS, A., FEBRERO, M. and FRAIMAN, R. (2002). Linear functional regression: The case of fixed design and functional response. *Canad. J. Statist.* **30** 285–300. MR1926066

[8] DAROLLES, S., FLORENS, J.-P. and RENAULT, E. (2002). Nonparametric instrumental regression. Working paper, GREMAQ, Univ. Social Science, Toulouse.

[9] DELAIGLE, A. and GIJBELS, I. (2002). Estimation of integrated squared density derivatives from a contaminated sample. *J. R. Stat. Soc. Ser. B Stat. Methodol.* **64** 869–886. MR1979392

[10] DONOHO, D. L. (1995). Nonlinear solution of linear inverse problems by wavelet-vaguelette decomposition. *Appl. Comput. Harmon. Anal.* **2** 101–126. MR1325535

[11] EFROMOVICH, S. and KOLTCHINSKII, V. (2001). On inverse problems with unknown operators. *IEEE Trans. Inform. Theory* **47** 2876–2894. MR1872847

[12] FAN, J. (1991). On the optimal rates of convergence for nonparametric deconvolution problems. *Ann. Statist.* **19** 1257–1272. MR1126324

[13] FAN, J. (1993). Adaptively local one-dimensional subproblems with application to a deconvolution problem. *Ann. Statist.* **21** 600–610. MR1232507

[14] FERRATY, F. and VIEU, P. (2000). Dimension fractale et estimation de la régression dans des espaces vectoriels semi-normés. *C. R. Acad. Sci. Paris Sér. I Math.* **330** 139–142. MR1745172

[15] HALL, P. and HOROWITZ, J. L. (2005). Nonparametric methods for inference in the presence of instrumental variables. *Ann. Statist.* **33** 2904–2929. MR2253107

[16] HOROWITZ, J. L. and MARKATOU, M. (1996). Semiparametric estimation of regression models for panel data. *Rev. Econom. Stud.* **63** 145–168. MR1372250

[17] JOHNSTONE, I. M. (1999). Wavelet shrinkage for correlated data and inverse problems: Adaptivity results. *Statist. Sinica* **9** 51–83. MR1678881

[18] KUNDUR, D. and HATZINAKOS, D. (1998). A novel blind deconvolution scheme for image restoration using recursive filtering. *IEEE Trans. Signal Processing* **46** 375–390.

[19] LI, T. and HSIAO, C. (2004). Robust estimation of generalized linear models with measurement errors. *J. Econometrics* **118** 51–65. MR2030966

[20] NATTERER, F. (1984). Error bounds for Tikhonov regularization in Hilbert scales. *Applicable Anal.* **18** 29–37. MR0762862

[21] NEUBAUER, A. (1989). Tikhonov regularisation for nonlinear ill-posed problems: Optimal convergence rates and finite-dimensional approximation. *Inverse Problems* **5** 541–557. MR1009038

[22] NEWEY, W. K. and POWELL, J. L. (2003). Instrumental variable estimation of nonparametric models. *Econometrica* **71** 1565–1578. MR2000257

[23] NEWEY, W. K., POWELL, J. L. and VELLA, F. (1999). Nonparametric estimation of triangular simultaneous equations models. *Econometrica* **67** 565–603. MR1685723

[24] RAMSAY, J. O. and SILVERMAN, B. W. (1997). *Functional Data Analysis*. Springer, New York.

[25] RAMSAY, J. O. and SILVERMAN, B. W. (2002). *Applied Functional Data Analysis: Methods and Case Studies*. Springer, New York. MR1910407
[26] SHENK, J. S. and WESTERHAUS, M. O. (1991). Population definition, sample selection, and calibration procedures for near infrared reflectance spectroscopy. *Crop Science* **31** 469–474.
[27] STANISWALIS, J. G. and LEE, J. J. (1998). Nonparametric regression analysis of longitudinal data. *J. Amer. Statist. Assoc.* **93** 1403–1418. MR1666636
[28] STEFANSKI, L. and CARROLL, R. J. (1990). Deconvoluting kernel density estimators. *Statistics* **21** 169–184. MR1054861
[29] TIKHONOV, A. N. (1963). On the solution of incorrectly put problems and the regularisation method. In *Outlines Joint Sympos. Partial Differential Equations (Novosibirsk, 1963)* 261–265. Acad. Sci. USSR Siberian Branch, Moscow. MR0211218
[30] VAN ROOIJ, A. C. M. and RUYMGAART, F. H. (1996). Asymptotic minimax rates for abstract linear estimators. *J. Statist. Plann. Inference* **53** 389–402. MR1407650
[31] WESLEY, I. J., UTHAYAKUMARAN, S., ANDERSSEN, R. S., CORNISH, G. B., BEKES, F., OSBORNE, B. G. and SKERRITT, J. H. (1999). A curve-fitting approach to the near infrared reflectance measurement of wheat flour proteins which influence dough quality. *J. Near Infrared Spectroscopy* **7** 229–240.
[32] YANG, Y., GALATSANOS, N. P. and STARK, H. (1994). Projection-based blind deconvolution. *J. Optical Soc. Amer. Ser. A* **11**, 2401–2409.

DEPARTMENT OF MATHEMATICS AND STATISTICS
THE UNIVERSITY OF MELBOURNE
MELBOURNE, VICTORIA 3010
AUSTRALIA
E-MAIL: halpstat@ms.unimelb.edu.au

DEPARTMENT OF ECONOMICS
ANDERSON HALL
NORTHWESTERN UNIVERSITY
2001 SHERIDAN ROAD
EVANSTON, ILLINOIS 60208-2600
USA
E-MAIL: joel-horowitz@northwestern.edu